
\baselineskip=14pt
\parskip=10pt

\font\eightrm=cmr8 

\magnification=\magstephalf

\def\1{{\overline{1}}}
\def\2{{\overline{2}}}
\parindent=0pt
\overfullrule=0in

\def\frac#1#2{{#1 \over #2}}
\centerline
{\bf 
Automating John P. D'Angelo's method to study Complete Polynomial Sequences
}
\bigskip
\centerline
{\it Shalosh B. EKHAD and Doron ZEILBERGER}
\bigskip
\qquad 
{\it In memory of Ron Graham\footnote{$^1$}
{\eightrm Ron Graham was not only a great combinatorialist, but also the ``coolest" (i.e. socially skilled) mathematican I have ever known. I confess
that for a long time, while I did admire his many contributions, I felt that he got more recognition than many other `nerdier' (i.e. socially awkward) mathematicians
who are even better than him, and, to my shame, I used him as an example of an `over-rated mathematican'. But this was wrong of me! In particular the
paper [G] that he wrote two years after his PhD, testifies to his great depth!}
 (1935-2020) the coolest mathematician  that I have ever known}

\bigskip

{\bf Abstract}: In a recent intriguing article by  complex geometer (and several complex variabler) John P. D'Angelo
there was a surprising application of elementary (but very deep!) number theory to complex geometry.
What D'Angelo needed was the largest integer {\bf not} representable as sum of distinct values of triangular numbers
(with analogous questions about other polynomial sequences). He demonstrated, in terms of a few lucid examples, how
to determine these numbers, and then rigorously prove that they are indeed correct. In this short note
we describe a Maple package that implements these ideas enabling, at least in principle, but often also
in practice, to determine these numbers {\bf fully automatically}. In fact, we show that
the very same ideas can be turned into an algorithm that inputs an arbitrary integer-generating polynomial,
and positive integers $a$ and $C$,  and outputs the smallest integer not representable as a sum
of distinct values of that polynomial with argument $\geq a$ in at least $C$ different ways. 
(The cases that D'Angelo needed were $a=0$ or $a=1$, and $C=1$).

{\bf Maple Package}

This article is accompanied by a Maple package {\tt JPDA.txt}, available from

{\tt https://sites.math.rutgers.edu/\~{}zeilberg/tokhniot/JPDA.txt} \quad .

The web-page of this article,

{\tt https://sites.math.rutgers.edu/\~{}zeilberg/mamarim/mamarimhtml/jpda.html} \quad ,

contains numerous input and output files, some ot them referred to in this paper.

{\bf Introduction}

In a recent fascinating article [D], John D'Angelo discovered a surprising connection between certain dimensions in complex geometry and
a problem in number theory tackled many years ago by such giants as Klaus Roth, George Szekeres, and Ron Graham [G].

The initial problem was, given a polynomial $p(j)$ , in $j$, that always yields positive integers for $j \geq 0$, to decide
whether every sufficiently large integer $n$ can be written as a sum of {\bf distinct} values of $p(j)$. In other
words, there exists an integer, {\it the largest exceptional value}, denoted (in [G]) by $\lambda(p)$, such that
$\lambda(p)$ is {\bf not representable} as a sum of distinct values, but every integer $n>\lambda(p)$ is representable in such a way.
Such polynomial sequences are called (in [D]) {\it Complete polynomial sequences}, and (in [G]) {\it Complete sequences of polynomial values}.

At the end of [G], Ron Graham commented

{\it ``$\dots$ Even for the simplest polynomials $p(j)$, the exact determination of $\lambda(p)$ is not easy''}

Already in 1948 (probably without computers) R. Sprague proved (see [G]) that
$$
\lambda(j^2)=128 \quad .
$$

In 1964 (most probably using computers), Graham proved that
$$
\lambda(j^3)=12758 \quad,
$$
$$
\lambda(j^4) > 2400000 \quad .
$$

How much do we know, {\it exactly}, today, with Moore's exponential growth of computer power? Not much more!

If you google (make sure to have the quotation marks) ``128, 12758'', you are immediately lead to OEIS sequence {\bf A001661} ({\tt https://oeis.org/A001661}),
that tells you that, in addition to the above values, known back in 1964, we currently know that
$$
\lambda(j^4) =5134240 \quad , \quad
\lambda(j^5) = 67898771 \quad , \quad
\lambda(j^6) = 11146309947 \quad , \quad
\lambda(j^7) = 766834015734 \quad .
$$
Surprisingly, except for the last number, these numbers already appeared in the 1995 edition [SP] {\bf M5393}, but only the first three numbers are listed in the first edition [S] (published 1973), where it is
sequence {\bf 2342}.

Ron Graham [G] gave a beautiful necessary and sufficient condition, easily checked (it is implemented in procedure {\tt RGc} in the Maple package {\tt JPDA.txt} accompanying our article)
for the polynomial $p(j)$ to be complete. In other words for the {\bf existence} of $\lambda(p)$. His proof entails an explicit upper bound, but as he commented himself,
it is {\it useless} in practice.

In [D], John D'Angelo gave a nice algorithm for rigorously determining the exact value of $\lambda(p)$ for {\it any given} polynomial $p(j)$, but rather than
formulating it formally, he illustrated it in terms of a few illuminating examples. In this short note, we will `reverse engineer' his argument, and
phrase it formally, and more importantly, {\bf implement} it, so that anyone interested in finding the exact value of $\lambda(p)$ (and hence the
exact dimension of some problems in complex geometry, see [D]) can do it by herself (or himself, or itself). Of course, it can only go so far,
and we  doubt that $\lambda(j^{20})$ will {\it ever} be known by human-kind (or machine-kind, for that matter).

In fact the algorithm easily generalizes to the more general scenario.

{\bf Inputs}: 

$\bullet$ A polynomial $p(j)$ that, thanks to Graham's condition (Theorem [4] in [G], Theorem 3.1 in [D]), is known to be {\it complete} \quad ;

$\bullet$ a non-negative integer $j_0$ \quad ;

$\bullet$ a positive integer $C$ \quad .

{\bf Output}: The largest integer, let's call it $\lambda_{j_0,C}(p)$, that is {\bf not} representable as sum of
{\bf distinct} values of $p(j)$ with $j \geq j_0$ in {\bf at least} $C$ different ways.

D'Angelo was particularly interested in $\lambda_{0,1}({{j+k} \choose {k}})$, $\lambda_{1,1}( {{j+k}\choose{k}})$ for $k=2,3,\dots$,
and gave the exact values for $k=2$ and $k=3$. He also asked about the ratio  of the latter quantities as $k$ goes to infinity.

In [D], the following exact values are given
$$
\lambda_{0,1}({{j+2} \choose {2}}) \, = \, 33 \quad, \quad \lambda_{1,1}({{j+2} \choose {2}}) \, = \, 50 \quad, \quad
$$
$$
\lambda_{0,1}({{j+3}\choose{3}}) \, = \, 558 \quad, \quad \lambda_{1,1}({{j+3} \choose {3}}) \, = \, 897 \quad .
$$
We further have (using D'Angelo's method, to be described shortly)
$$
\lambda_{0,1}({{j+4} \choose {4}}) \, = \, 12659 \quad, \quad \lambda_{1,1}({{j+4} \choose {4}}) \, = \, 23319 \quad, \quad
$$
$$
\lambda_{0,1}({{j+5} \choose {5}}) \, = \, 120838 \quad, \quad \lambda_{1,1}({{j+5} \choose {5}}) \, = \, 291217 \quad.
$$

The sequence $33,558,12659,120838,\dots$, that may be viewed as the ${ {j+k} \choose{k} }$ analog of OEIS sequence {\bf A001661} (that concerns $j^k$), (and in some sense is a more natural one, since in discrete math
the natural basis of polynomials is the former, rather than the latter) was not (viewed Nov. 3, 2021) in the OEIS. We hope that some kind soul would enter it.
We also recommend the companion sequence (where the case $j=0$, i.e. $1$, is not allowed to participate), $50,897,233319,291217,\dots$, that is also of interest to John D'Angelo, since it has complex geometry implications.

\vfill\eject

{\bf Formalizing John D'Angelo's algorithm}

{\bf Step 1}: Fix a positive integer $K$ for exploration, and using Maple, or any other computer algebra system, find the beginning 
up to the power $q^K$ of the generating function
$$
\prod_{j=j_0}^{\infty} (1+ q^{p(j)}) \quad .
$$
By the usual {\it generatingfunctionolgy} argument, the coefficient of $q^n$ is the exact number of ways of representing $n$ as a distinct sum of values of $p(j)$ with $j \geq j_0$.

Find the largest $n \leq K$ whose coefficient is $<C$. We conjecture that this is also the largest integer $n< \infty$ with that property, but so far it is only a conjecture.
Let's call this value $N_0$.

in order to prove our conjecture that $\lambda_{j_0,C} (p) =N_0$,  {\it rigorously}, we formulate the following statement, regarding a positive integer $y$:

$S(y):=$ $y$ is representable as a sum of distinct values of $p(j)$ for $j \geq j_0$ in at least $C$ different ways.

We  believe (based on the initial numerical exploration) that $S(y)$ is true for all $y > N_0$, but we would like to prove it rigorously.
By mathematical induction, if we can prove that $S(y')$ for $N_0<y'<y$ implies $S(y)$, we will be done. In all inductive proofs we need
a {\it base case}, or more often {\it base cases}. Let's leave them open for now. In other words, in analogy to {\it undetermined coefficients}, we have
{\it undetermined base cases for the inductive argument}. So far we know (from the initial exploration, that $S(y)$ is true for $N_0<y \leq K$, but
perhaps that does not suffice. Let's see.

Given such a positive integer $y$, let $N$ be the unique integer such that
$$
p(N)<y \leq p(N+1) \quad.
$$

D'Angelo's ingeniously writes
$$
y= (y-p(N-3))+ p(N-3) \quad .
$$

Now, {\it wishful thinking}, {\bf if} we knew the following two facts
$$
y-p(N-3) > N_0 \quad,
$$
$$
y-p(N-3)<p(N-3) \quad ,
$$
then things would follow by induction. By the inductive hypothesis, $y-p(N-3)$ is representable as a sum of values of $p(j)$ in at least $C$ different ways,
and for each of them we can add $p(N-3)$ to get such a representation for $y$. The second condition guarantees that appending $p(N-3)$ will not
violate uniqueness (since there is no way that $p(N-3)$ can be a part of a representation of $y-p(N-3)$).

But how to guarantee these conditions?

We have
$$
y-p(N-3)=(y-p(N))+(p(N)-p(N-3))> p(N)-p(N-3)
$$
So if we knew that $p(N)-p(N-3) \geq N_0$, then we would be in good shape.

In order to satisfy the second inequality, note that

$$
y-p(N-3) \leq  p(N+1)-p(N-3)<p(N-3) \quad,
$$

provided $p(N+1)<2p(N-3)$. So we also need

$$
2p(N-3) -p(N+1)>0 \quad.
$$

Both polynomials $p(N)-p(N-3)$  and $2p(N-3) -p(N+1)$ have positive leading coefficient hence there is a {\bf cutoff} where  they both start being positive. Let that cutoff be $N_1$. Then the above
inductive argument  is valid for $y>p(N_1)$. This leaves, as {\bf base cases} for the induction argument, all the integers from $N_0$ to $p(N_1)$. If $p(N_1)$ is less than
the $K$ that we used above for exploration, we are already done. Otherwise, let the computer check that the coefficients of the generating function $\prod_{j=j_0}^{\infty} (1+ q^{p(j)})$
between $N_0$ and $p(N_1)$ are all $\geq C$.

The web-page of this paper contains numerous sample input and output files. In particular procedure {\tt CutoffV(p,j,st,C,K)}
in the Maple package not only finds the cutoffs (and checks them), but spells out the above argument for each given case. See for example:

{\tt https://sites.math.rutgers.edu/\~{}zeilberg/tokhniot/oJPDA1.txt} \quad .

{\bf References}

[D] John P. D'Angelo, {\it Symmetries, Rational Sphere Maps, and Complete Polynomial Sequences}, preprint, submitted for publication. \hfill\break
{\tt https://faculty.math.illinois.edu/\~{}jpda/jpd-symmetries-2.pdf} \quad .

[G] Ron L. Graham, {\it Complete sequences of polynomial values}, Duke Math. J. {\bf 31} (1964), 275-285.
{\tt https://sites.math.rutgers.edu/\~{}zeilberg/akherim/graham1964.pdf} \quad.

[S] N. J. A. Sloane, {\it ``A Handbook of Integer Sequences''}, Academic Press, 1973

[SP] N. J. A. Sloane and Simon Plouffe, {\it ``The Encyclopedia of Integer Sequences''}, Academic Press, 1995.

\bigskip
\hrule
\bigskip
Shalosh B. Ekhad and Doron Zeilberger, Department of Mathematics, Rutgers University (New Brunswick), Hill Center-Busch Campus, 110 Frelinghuysen
Rd., Piscataway, NJ 08854-8019, USA. \hfill\break
Email: {\tt [ShaloshBEkhad, DoronZeil] at gmail dot com}   \quad .

{\bf Exclusively published in the Personal Journal of Shalosh B. Ekhad and Doron Zeilberger and arxiv.org}

Written: {\bf  Nov. 3, 2021}.

\end